# Gender differences in grant peer review:

# A meta-analysis

Lutz Bornmann,* Ruediger Mutz, Hans-Dieter Daniel

Addresses for all authors

First author:
Lutz Bornmann*
Professorship for Social Psychology and Research on Higher Education
ETH Zurich
Zaehringerstr. 24
CH-8092 Zurich
Phone: +49 44 632 4825
Fax: +49 44 632 1283
E-mail: bornmann@gess.ethz.ch

Second author:
Ruediger Mutz
Professorship for Social Psychology and Research on Higher Education
ETH Zurich

Third author:
Hans-Dieter Daniel
Evaluation Office
University of Zurich

Professorship for Social Psychology and Research on Higher Education
ETH Zurich

**Abstract**

Narrative reviews of peer review research have concluded that there is negligible evidence of gender bias in the awarding of grants based on peer review. Here, we report the findings of a meta-analysis of 21 studies providing, to the contrary, evidence of robust gender differences in grant award procedures. Even though the estimates of the gender effect vary substantially from study to study, the model estimation shows that all in all, among grant applicants men have statistically significant greater odds of receiving grants than women by about 7%.

Keywords

meta-analysis, grant peer review, gender bias, gender effect, gender differences

# Introduction

Robert K. Merton (1973) describes the "norm of universalism" in science whereby all scientists should be evaluated purely on the quality of their work and thus have the same chance of success. Merton (1973) regarded meritocracy in science as essential for the advancement of scientific knowledge. If functionally irrelevant factors skew the results when scientific contributions are evaluated, it will take longer for nature to reveal its laws (Hull, 1990). A functionally irrelevant factor is one that does not affect the ability of an individual to perform a particular functional role or status. For the progress of science, a grant applicant's gender is functionally irrelevant. If gender is used, either explicitly or implicitly, in the evaluation of scientific work and affects the way in which grant decisions are made, the principles of universalism and objectivity are being abridged (Cole, 1992). Potential disadvantages for women compared with men have been and remain a major source of controversy in society, and scientists have been heavily preoccupied with the issue in the context of the operation of the grant awards system (Wennerås & Wold, 2000).

Even if most scientific institutions have affirmative action officers and procedures designed to ensure that women have an equal chance, a comprehensive review of the literature on gender differences in the careers of academic scientists by the National Science Foundation (NSF, Arlington, VA, USA) comes to conclusion: "Taken as a whole, the body of literature we reviewed provides evidence that women in academic careers are disadvantaged compared with men in similar careers. Women faculty earn less, are promoted less frequently to senior academic ranks, and publish less frequently than their male counterparts" (National Science Foundation, 2003, p. 1). In the award of grants based on peer review, on the other hand, the gender of the applicant appears to play a negligible part. In 2003 a literature review of peer review research concluded: "The most frequent criticism made by scientists about the day to day operation of peer review is that of … gender bias. We suggest that this criticism is generally unfounded" (Wood & Wessely, 2003, p. 36). In a Cochrane Methodology Review the research findings are summarized in similar terms: "Descriptive evidence of gender bias was provided by a study at the Swedish Medical Research Council, but a number of other studies carried out in similar contexts found no evidence of it" (Demicheli & Pietrantonj, 2004).

Both of the above-mentioned reviews of the research on peer review finding little evidence of gender bias describe the existing literature using the narrative technique, without attempting quantitative synthesis of study results. As from the viewpoint of quantitative social scientists narrative reviews are not very precise in their descriptions of study results (Shadish, Cook, & Campbell 2002), quantitative techniques should be used as well as narrative techniques. The term "meta-analysis" (Glass, 1976) refers to a statistical approach that combines evidence from different studies to obtain an overall estimate of treatment effects (Shadish, Cook, & Campbell, 2002). Taking the quantitative approach, meta-analysis allows generalized statements on the strength of the effects, regardless of the specificity of individual studies (Matt & Navarro, 1997). Undertaking a meta-analysis presupposes that each study is similarly designed with regard to certain properties (e.g., methods, sampling). Even though the studies that investigated the influence of the applicant's gender on decisions in grant peer review are quite heterogeneous in these properties, most of them are very similar in reporting the following proportions: number of men and number of women among both approved and rejected applicants for grants. Using these frequency data, it is possible to evaluate the empirical studies meta-analytically (Skrondal & Rabe-Hesketh, 2004).

Meta-analysis (Skrondal & Rabe-Hesketh, 2004) is recommended for a quantitative review of experimental studies in which subjects were randomly assigned to a treatment or control group. However, peer review is not a 'medication' that can be tested in a randomized clinical trial (Marusic, 2005). For the meta-analysis we were therefore restricted to non-experimental studies that investigated peer review in a "natural setting." The restriction to the "natural setting" makes it very difficult to establish unambiguously whether work from a particular group of scientists has a higher approval rate due to biases[1] in the decision-making procedure, or if greater success in the procedure is simply a consequence of the scientific merit of the corresponding group of applications. As with non-experimental studies it is hardly possible to demonstrate gender bias in the decision-making process of grant peer review, our meta-analysis addressed the following question: Do the studies show an overall gender *effect* in the decisions in grant peer review?

## Methods

*Locating studies*

The literature research was conducted at the end of 2005. In a first step we located some studies that investigated the influence of the applicant's gender on decisions in grant peer review, using the reference lists provided by narrative reviews of research on grant peer review (Demicheli & Pietrantonj, 2004; Wood & Wessely, 2003) and using tables of contents of certain journals (e.g., *Nature, Research Policy, Scientometrics*). We conducted a search of both publications (journal articles, monographs, collected works, etc.) and grey literature (Internet documents, institutional reports, case reports, etc.) to avoid a bias that follows from the difficulties of publishing nonconforming studies (Eagly, 2005; Greenwald, 1975).

---

[1] "Bias is any feature of an evaluator's cognitive or attitudinal mind-set that could interfere with an objective evaluation" (Shatz, 2004, p. 36).

In a second step, in order to obtain keywords for searching computerized databases, we prepared a bibliogram (White, 2005) for the studies located in the first step. The bibliogram ranks by frequency the words included in the abstracts of the studies located. Words at the top of the ranking list (e.g., *peer review*, *gender*, *sex*, *women*, *science,* and *bias*) were used for searches in computerized literature databases (e.g., *Web of Science*, *IngentaConnect*, *PubMed*, *Sociological Abstracts*, *ProQuest Digital Dissertations*, *PsycINFO, ERIC*) and Internet search engines (e.g., *Google, Ask Jeeves*).

In the final step of our literature search, we located all of the citing publications for a series of articles (found in the first and second step) for which there are a fairly large number of citations in *Web of Science*. In addition, we sent inquiries to authors of narrative reviews on peer review research asking them to supplement our publications list with further studies. If the data in the studies that we researched was incomplete for our purposes, we tried to contact the authors of the studies. But it is often difficult to trace the authors of older studies, and they are often hazy even about relatively recent work (Shadish, Cook, & Campbell, 2002).

We encountered some difficulty in finding studies, because grant peer review – although "the key institution of the scientific culture" (Ziman, 2000, p. 246) – frequently remains a secret activity (Tight, 2003). According to Daryl Chubin, director of the Center for Advancing Science & Engineering Capacity (American Association for the Advancement of Science, Washington, DC), the reason peer review is seldom studied empirically is that reviews are secured with assurance of confidentiality. To link reviewer and review information would break that confidence. Where a systematic review of peer review is undertaken, the findings are often not published. If the findings are published, their reception can be hampered if the relevant documents are available only in the respective native language.

*Studies for the meta-analysis*

We were able to include in our meta-analysis data on proportions of women and men for 66 different peer review procedures from 21 studies. As the studies give these proportions for various funding institutions, funding programs, and different cohorts of applicants, we had data for three different procedures on average per study. Thus, one procedure is characterized by a single set of proportions: women and men among applicants as well as women and men among approved applicants. For example, Ackers et al. (2000) report the participation of women researchers in the TMR Marie Curie Fellowship Programme. The authors present the proportions of women and men for seven peer review procedures (that is, chemistry, engineering, mathematics, earth sciences, economics, physics, and life sciences). Table 1 (part 1) and Table 2 (part 2) provide the information collected by our research team on the 66 peer review procedures from the studies. One researcher on our team extracted this information from the publications. The information was then validated by a second researcher on our team, and disagreements were resolved by consensus. The bibliographic data for the 21 studies can be found in Table 3.

The studies included in the meta-analysis were published between 1987 and 2005 as article ($n = 10$), grey literature ($n = 7$), letter ($n = 2$), book ($n = 1$), or commentary ($n = 1$) by men scientists ($n = 9$), women scientists ($n = 6$), or by both men and women scientists ($n = 5$) (the authors of one study are not known). In one study the proportion of women and men is based on survey data; the other studies worked with archive data. The data evaluated in the studies refers to periods between 1979 and 2004. Due to the fact that most of the studies used were published as peer reviewed publications or as reports of renowned funding institutions their quality is sufficiently guarantied.

The data on the 66 peer review procedures in the 21 studies refer either to peer review for the assessment of research applications ($n = 40$) or peer review for the assessment of applications for post-graduate fellowships ($n = 26$). The studies covered 13 findings on the proportion of women and men for a research funding institution in Australia (e.g., Australian Research Council, Canberra); 15 findings refer to the peer review of a funding institution in North America (e.g., NSF) and 38 that of a European institution (e.g., EMBO). Twenty-nine sets of data in the studies refer to peer review in the life sciences field; for the exact sciences there are 11 sets and for the social sciences and humanities 13 sets of data. Twelve proportion findings concern more than one specialist discipline (in one study there is no classification of discipline).

*Statistical methods*

For the calculation of the meta-analysis we considered as dependent variable estimations of the odds ratio, the odds for being approved among women applicants divided by the odds of being approved among men applicants. For an individual peer review procedure $j$ (e.g., at the NSF) that we included in the meta-analysis, this odds ratio can be estimated as

$$o_j = \frac{d_{1j}/(n_{1j} - d_{1j})}{d_{0j}/(n_{0j} - d_{0j})}, \tag{1}$$

where $d_{1j}$ and $n_{1j}$ are the number of women among approved and all applicants, respectively, whereas $d_{0j}$ and $n_{0j}$ are the number of men among approved and all applicants, respectively.

The approach of our meta-analysis is to analyze the estimated gender effect at the level of different peer review procedures that were investigated in the different studies. As we assume that the true gender effect varies between peer review procedures and between studies, we used for the meta-analysis *generalized linear mixed models* that explicitly allow for these variations in a multilevel framework. As the estimates of effect sizes violate the normality assumption, we analyzed log-odds ratios instead of odds ratios (Skrondal & Rabe-Hesketh, 2004).

In our data there is a three-level structure with applicants nested within procedures and procedures nested within studies. The predominant approach in meat-analysis is to analyze the procedure or study effects instead of the original applicant-level data (Skrondal & Rabe-Hesketh, 2004, p. 301). For the data, where *j* indicates procedures, and *k* studies, we consider the generalized linear mixed model of the log-odds ratio ($o_j$), as defined in Eq. (1).

$$\text{logit}(o_j) = \beta_0 + \beta_1\, x_{jk} + \zeta_{0jk} + \zeta_{1jk}\, x_{jk} + \zeta_{0k} + \zeta_{1k}\, x_{jk} + \varepsilon_j \qquad \varepsilon_j \sim \text{N}(0, \theta_j) \qquad (2)$$

where

$$x_{jk} = \begin{cases} +\,0.5\ \textit{for women applicants} \\ -\,0.5\ \textit{for men applicants} \end{cases},$$

and

$$(\zeta_{0j},\, \zeta_{1j}) \sim \text{N}(\mathbf{0},\, \mathbf{\Psi}^{(2)})$$

$$(\zeta_{0k},\, \zeta_{1k}) \sim \text{N}(\mathbf{0},\, \mathbf{\Psi}^{(3)})$$

Here $\beta_0$ is a fixed intercept and $\zeta_{0jk},\, \zeta_{0k}$ are random intercepts; $\beta_1$ is a fixed slope, and $\zeta_{1jk},\, \zeta_{1k}$ are random slopes of $x_{jk}$. $\beta_1$ represents the log odds ratio of interest (gender effect), whereas $\beta_1 + \zeta_{1jk} + \zeta_{1k}$ represents the estimated 'true' log odds ratio of the peer review procedure *j* examined in study *k*. The within-study deviations $\sqrt{\theta_j}$ are simply set equal to the standard errors of the log-odds ratios, calculated out of the proportions given in Eq. 1 (Woolf's method). The sum of the terms $\zeta_{0jk} + \zeta_{1jk} x_{jk} + \zeta_{0k} + \zeta_{1k} x_{jk}$ (the random part) can be thought of as a total residual $\xi_{jk}$. The variance of this total residual var($\xi_{jk}$) is composed of six variance and covariance components: the variances of the intercepts between procedure $\zeta_{0jk}$ and between studies $\zeta_{0k}$, the variances of slopes between procedures $\zeta_{1jk} x_{jk}$ and between studies $\zeta_{1k} x_{jk}$, and the two covariances between intercepts and slopes on the two levels: cov($\zeta_{0jk}$, $\zeta_{1jk} x_{jk}$), cov($\zeta_{0k}$, $\zeta_{1k} x_{jk}$). Thus, two matrices of variance/covariance components can be generated: one for the random effects of the procedures ($\mathbf{\Psi}^{(2)}$) and one for the random effects of the studies ($\mathbf{\Psi}^{(3)}$). The generalized linear mixed model was estimated by maximum likelihood using adaptive quadrature. For the analysis the gllamm-procedure (Rabe-Hesketh, Skrondal, & Pickles, 2004) of the software STATA 9.0 (StataCorp., 2005) was used.

## Results

Table 4 shows the maximum likelihood estimates for the generalized linear mixed model that examined gender effect in grant peer review. The fixed part represents the overall regression on gender (men = -0.5; women = +0.5); the random part represents the variability of regressions between procedures and between studies, as well as the heterogeneity of effects (expressed as variances and covariances of intercepts and slopes). The estimate of $\beta_1$ (-0.07) in Table 4 (fixed part) gives, over all 66 peer review procedures analyzed in the 21 studies, the effect of the applicant's gender on the approval (i.e., funding) of grant applications. There appears to be strong evidence that the odds of being approved for a grant are smaller for women than for men scientists, with an estimated odds ratio of exp($\beta_1$ -0.5 – $\beta_1$ 0.5) = exp(0.07) = 1.07. In other words, in grant peer review procedures, men have on average statistically significantly greater odds of approval than women applying for grants by about 7%. Due to the high number of different peer review procedures ($n$=66, in terms of funding institutions, funding programs, fields of promotion, etc.) and reviewed applications ($n$=353725) in our meta-analysis, the greater odds of approval for men compared to odds of approval for women can be interpreted as a *robust* average gender effect in grant peer review (Shadish, Cook, & Campbell, 2002).

The findings concerning the random part of the model estimation (see Table 4) suggest that there is a statistically significant variability in the gender effect (preference for men over women) both at the level of the peer review procedures (0.04) and at the level of the studies (0.07). Accordingly, even if a statistically significant gender effect can be seen over all studies, the size of the effect in the individual peer review procedures varies widely. The $z$-value of the corresponding variance components var($\zeta_{0jk}$), var($\zeta_{1jk}$) for procedures and var($\zeta_{0k}$), var($\zeta_{1k}$) for studies are significant. Figure 1 shows the variability of the gender effect between peer review procedures using "Empirical Bayes" estimates as effect sizes and credible intervals for each peer review procedure (the figure shows both the names of the studies and the research funding institutions (abbreviated) to which the peer review predictions refer). The proponents of Bayesian statistics (Carlin & Louis, 2000) use credible intervals instead of confidence intervals in order to handle such complex structures (Skrondal & Rabe-Hesketh, 2004). For a 95% credible interval, the posterior probability that the gender effect parameter for a special procedure lies in the interval is 95%. In Figure 1 negative mean log odds ratios point to a preference for men, and positive mean log odds ratios to a preference for women in the peer review. The estimated individual mean log odds ratios vary between -0.26 (Ackers, 2000a; Marie Curie) and 0.20 (Ackers, 2000d; Marie Curie). In other words, the gender effects (odds ratios) vary over the studies within the range of 22.1% (exp(-0.26) =

1.221) in favor of men and 22.9% (exp(-0.20) = 0.771) in favor of women applying for grants (however, Figure 1 pointed clearly out that there are only a few procedures with an effect in favor of women).

## Discussion

Conventionally, peer review is regarded as a sure guarantee of good science. It reassures us about the quality of scientific work and that the taxpayer's money spent on science is well spent (Biagioli, 2002). The findings of our quantitative review suggest, in contrast to recently published narrative reviews of the research on grant peer review (Demicheli & Pietrantonj, 2004; Wood & Wessely, 2003), that there are overall robust gender difference in grant peer review procedures. Although the parameter estimations of the gender effect vary substantially for the individual procedures, the model estimation shows that among applicants for grants, men have greater odds of approval than women by about 7%. In other words, in grant peer review the odds for approval of men's grant applications are about 15:14 (for women the odds are reversed). This finding tallies with the conclusion of the National Science Foundation (2003) on gender differences in the careers of academic scientists.

What importance can we attribute to this finding on the use of peer review for selection of applicants to receive grants? Assume for a moment that worldwide in a certain time period decisions are made to approve and reject 100,000 grant applications (50,000 submitted by women and 50,000 submitted by men scientists). Half of these applications are approved, and half are rejected. Based on the results of our meta-analysis, we can expect an approval rate for grant applications submitted by men scientists of 52% (26,000 approvals). For grant applications submitted by women, we can expect an approval rate of 48% (24,000 approvals). This makes a difference of 2,000 approvals that is due to the gender of the applicants.

However, the cause of this discrepancy is unknown. We tried to include characteristics of the peer review procedures (see Table 2) into the model estimation to get some hints for potential causes. For instance, the proportions of women and men for the procedures are highly aggregated over fields of study or disciplines, respectively. One could easily imagine that women are more likely to be represented in the softer sciences (such as social sciences) in which the overall success rates are lower than those in the harder sciences (such as chemistry). Even if there are no gender differences within fields of study, aggregation over fields of study can create a strong effect ("Simpson's paradox"). Another characteristic of the peer review procedures that might be of interest for inclusion into the model estimation is the cohort of application. The data evaluated in the 21 studies refers to periods between 1979 and 2004. It is safe to assume that in this period of time considerable progress has been made to reduce gender bias in science and science funding.

Unfortunately, the inclusion of one or more of these characteristics into the calculation of the meta-analysis resulted in models that did not converge in the estimation process. This finding indicated that the model estimation became too complex by considering specific interaction effects or the included characteristics had no influence on the outcome, respectively.

In the literature on peer review a number of possible causes for gender discrepancies are suggested or ruled out by research findings. According to the study of Wennerås and Wold (1997) the differences are *not* the result of differences in scientific productivity before applying for a grant (see also Bornmann & Daniel, 2005). The findings of this well received study point out that we need other explanations. Instead, gender differences in grant peer review procedures could be attributed to fewer women principal investigators applying for grants. Gender bias – whether implicit or explicit – could come into play. According to Valien (1999) "one remedy for the salience of sex … is having adequate numbers of women and men in any group being evaluated. That will make any given woman's sex less distinctive" (p. 309). The explanation for gender differences could be also institutional; there are more men in higher-ranking positions than women (more men held senior rank positions than women), meaning fewer women are in decision-making positions.

Whatever the cause, there are ways to rule out gender bias – whether intentional or unintentional – from the grant-making process (for an overview about the wide range of ways, see Valian, 1999): One possible way is to mask applicants' gender (Bornmann, 2007). For example, it was found that "introducing a screen to obscure the gender of musicians auditioning for symphony orchestra positions increased the likelihood that a woman was selected by 30 to 60%" (Handelsman et al., 2005, p. 1191). In journal peer review, masking authors' gender has proved to be a satisfactory precaution against gender as a potential source of bias (Horrobin, 1982). Masking is not equally suitable for all document types of submissions, however. It is essentially impossible to pass valid judgment on a manuscript of the short-communication type without some personal knowledge of the author (Daniel, 1993, 2004; Ziman, 1968).

It is questionable as well whether the gender of the applicant in grant peer review should be masked. Apart from assessment of the grant applicant's proposed research, decisions on applications are also based on an assessment of the applicant's track records (Cox, Gleser, Perlman, Reid, & Roeder, 1993). According to Gerlind Wallon, program manager at the European Molecular Biology Organization (EMBO, Heidelberg, Germany), EMBO is planning an experiment in which they will totally gender-blind the committee in the research fellowship selection process (the results will be available by the end of 2007). As the example of the Pioneer Award for innovative research of the National Institutes of Health (NIH, Bethesda, MD, USA) shows, there are other effective measures against gender as a potential source of bias. While in 2004 there were no women among the nine winners (Mervis, 2004), in 2005 about half of the 13 scientists chosen for the Pioneer Award were women (Mervis, 2005). The improved chances for women scientists were said to be attributable to the fact that, for example, (1) women were especially encouraged to apply, (2) only self-nominations (rather than institutional submissions) were accepted, and (3) the NIH spent more time schooling its reviewers (Mervis, 2005).

Table 1. Characteristics of the peer review procedures in the studies (part 1)

| Study and research funding program/ institution | Document type of the study | Gender of the authors of the study | Number of submitted applications | Number of approved applications | Proportion of women among all applicants | Proportion of women among approved applicants | z value[1] |
|---|---|---|---|---|---|---|---|
| Ackers (2000a; Marie Curie) | grey literature | male and female | 1740 | 413 | 41 | 34 | 2,76 |
| Ackers (2000b; Marie Curie) | | | 1166 | 211 | 22 | 21 | 0,26 |
| Ackers (2000c; Marie Curie) | | | 1164 | 263 | 20 | 17 | 1,37 |
| Ackers (2000d; Marie Curie) | | | 758 | 159 | 33 | 40 | 1,54 |
| Ackers (2000e; Marie Curie) | | | 2028 | 409 | 45 | 38 | 2,45 |
| Ackers (2000f; Marie Curie) | | | 2833 | 574 | 21 | 20 | 0,51 |
| Ackers (2000g; Marie Curie) | | | 4302 | 870 | 47 | 44 | 1,80 |
| Allmendinger (2002a; DFG) | article | male and female | 85 | 61 | 15 | 13 | 0,37 |
| Allmendinger (2002b; DFG) | | | 90 | 58 | 9 | 9 | 0,06 |
| Allmendinger (2002c; DFG) | | | 105 | 69 | 8 | 9 | 0,25 |
| Allmendinger (2002d; DFG) | | | 110 | 61 | 15 | 13 | 0,26 |
| Allmendinger (2002e; DFG) | | | 103 | 47 | 11 | 9 | 0,43 |
| Allmendinger (2002f; DFG) | | | 137 | 75 | 32 | 27 | 0,84 |
| Allmendinger (2002g; DFG) | | | 131 | 75 | 11 | 7 | 1,19 |
| Bazeley (1998; ARC) | article | female | 400 | 93 | 14 | 12 | 0,58 |
| Bornmann (2005; BIF) | article | male | 2697 | 634 | 40 | 32 | 3,87 |
| Brouns (2000a; NWO/ KNAW) | article | female | 45 | 33 | 27 | 24 | 0,24 |
| Brouns (2000b; NWO/ KNAW) | | | 60 | 32 | 27 | 19 | 0,88 |
| Brouns (2000c; NWO/ KNAW) | | | 261 | 148 | 10 | 14 | 1,36 |
| Brouns (2000d; NWO/ KNAW) | | | 202 | 86 | 17 | 10 | 1,51 |
| Brouns (2000e; NWO/ KNAW) | | | 241 | 112 | 28 | 27 | 0,20 |
| Dexter (2002a; Wellcome Trust) | grey literature | male | 454 | 138 | 38 | 39 | 0,17 |
| Dexter (2002b; Wellcome Trust) | | | 1081 | 337 | 27 | 25 | 0,70 |
| Dexter (2002c; Wellcome Trust) | | | 87 | 34 | 20 | 12 | 1,12 |
| Emery (1992; NIH) | article | male and | 294 | 62 | 31 | 37 | 0,92 |

| Study and research funding program/ institution | Document type of the study | Gender of the authors of the study | Number of submitted applications | Number of approved applications | Proportion of women among all applicants | Proportion of women among approved applicants | z value[1] |
|---|---|---|---|---|---|---|---|
| | | female | | | | | |
| Friesen (1998a; MRC) | letter | male | 5535 | 1457 | 21 | 20 | 0,59 |
| Friesen (1998b; MRC) | | | 529 | 84 | 27 | 24 | 0,64 |
| Friesen (1998c; MRC) | | | 2219 | 333 | 39 | 33 | 1,92 |
| Goldsmith (2002a; NSF) | grey literature | female | 4288 | 524 | 34 | 28 | 2,82 |
| Goldsmith (2002b; NSF) | | | 8634 | 870 | 45 | 44 | 0,64 |
| Grant (1997a; Wellcome Trust) | letter | male | 1365 | 374 | 20 | 19 | 0,17 |
| Grant (1997b; Wellcome Trust) | | | 325 | 21 | 32 | 43 | 0,95 |
| Grant (1997c; MRC) | | | 784 | 212 | 21 | 23 | 0,56 |
| Grant (1997d; MRC) | | | 276 | 39 | 39 | 28 | 1,36 |
| Jayasinghe (2001; ARC) | article | male | 2981 | 635 | 15 | 15 | 0,02 |
| National Science Foundation (2005a) | grey literature | | 29928 | 9809 | 18 | 20 | 4,02 |
| National Science Foundation (2005b) | | | 28140 | 9261 | 20 | 21 | 1,92 |
| National Science Foundation (2005c) | | | 28337 | 9110 | 19 | 18 | 0,63 |
| National Science Foundation (2005d) | | | 29180 | 9727 | 19 | 20 | 2,48 |
| National Science Foundation (2005e) | | | 31349 | 9761 | 19 | 19 | 1,70 |
| National Science Foundation (2005f) | | | 34204 | 10215 | 20 | 20 | 0,22 |
| National Science Foundation (2005g) | | | 38573 | 10585 | 19 | 20 | 1,67 |
| National Science Foundation (2005h) | | | 41727 | 10041 | 20 | 21 | 1,99 |
| Over (1996; ARC) | article | male | 256 | 74 | 15 | 16 | 0,20 |
| Sigelman (1987; NSF) | article | male | 146 | 33 | 12 | 6 | 1,13 |
| Taplick (2005a; EMBO) | grey literature | male | 1070 | 581 | 54 | 53 | 0,46 |
| Taplick (2005b; EMBO) | | | 7703 | 1498 | 42 | 36 | 4,43 |
| Taplick (2005c; EMBO) | | | 1068 | 142 | 25 | 23 | 0,68 |
| Viner (2004a; EPSRC) | article | male | 5681 | 3378 | 3 | 4 | 1,07 |
| Viner (2004b; EPSRC) | | | 5431 | 3456 | 4 | 4 | 0,17 |

| Study and research funding program/ institution | Document type of the study | Gender of the authors of the study | Number of submitted applications | Number of approved applications | Proportion of women among all applicants | Proportion of women among approved applicants | z value[1] |
|---|---|---|---|---|---|---|---|
| Ward (1998; NHMRC) | article | male and female | 421 | 145 | 29 | 30 | 0,26 |
| Wellcome Trust (1997) | grey literature | male and female | 138 | 70 | 17 | 21 | 0,82 |
| Wenneras (1997; MRC) | commentary | female | 114 | 20 | 46 | 20 | 2,54 |
| Willems (2001a; DFG) | grey literature | female | 4687 | 2611 | 15 | 15 | 0,52 |
| Willems (2001b; DFG) | | | 5106 | 2813 | 16 | 16 | 0,75 |
| Willems (2001c; DFG) | | | 5048 | 2891 | 15 | 15 | 0,49 |
| Willems (2001d; DFG) | | | 5102 | 2871 | 17 | 17 | 0,60 |
| Wood (1997a; ARC) | book | female | 292 | 68 | 17 | 15 | 0,57 |
| Wood (1997b; ARC) | | | 337 | 73 | 17 | 12 | 0,99 |
| Wood (1997c; ARC) | | | 211 | 50 | 6 | 4 | 0,67 |
| Wood (1997d; ARC) | | | 191 | 46 | 7 | 2 | 1,80 |
| Wood (1997e; ARC) | | | 316 | 71 | 8 | 7 | 0,26 |
| Wood (1997f; ARC) | | | 363 | 82 | 29 | 30 | 0,33 |
| Wood (1997g; ARC) | | | 253 | 66 | 5 | 5 | 0,20 |
| Wood (1997h; ARC) | | | 480 | 124 | 25 | 26 | 0,09 |
| Wood (1997i; ARC) | | | 363 | 82 | 5 | 7 | 0,85 |

*Note.* [1] The *z* value provides a test of whether the proportion of women among approved applicants differs significantly from the proportion of women among all applicants. If *z* > 1.96, the two percentages are significantly different.

Table 2. Characteristics of the peer review procedures in the studies (part 2)

| Study and research funding program/ institution | Data source of the study | Cohort of application | Funding organization/ program that was examined in the study | Country location of the organization | Kind of promotion | Discipline or field |
|---|---|---|---|---|---|---|
| Ackers (2000a; Marie Curie) | archive data | 1994-1998 | Marie Curie | Europe | fellowship | chemistry |
| Ackers (2000b; Marie Curie) | | | | | | engineering |
| Ackers (2000c; Marie Curie) | | | | | | mathematics |
| Ackers (2000d; Marie Curie) | | | | | | earth sciences |
| Ackers (2000e; Marie Curie) | | | | | | economics |
| Ackers (2000f; Marie Curie) | | | | | | physics |
| Ackers (2000g; Marie Curie) | | | | | | life sciences |
| Allmendinger (2002a; DFG) | archive data | 1993 | German Research Foundation (DFG) | Germany | grant | sociology |
| Allmendinger (2002b; DFG) | | 1994 | | | | |
| Allmendinger (2002c; DFG) | | 1995 | | | | |
| Allmendinger (2002d; DFG) | | 1996 | | | | |
| Allmendinger (2002e; DFG) | | 1997 | | | | |
| Allmendinger (2002f; DFG) | | 1998 | | | | |
| Allmendinger (2002g; DFG) | | 1999 | | | | |
| Bazeley (1998; ARC) | archive data | 1995 | Australian Research Council (ARC) | Australia | grant | across disciplines |
| Bornmann (2005; BIF) | archive data | 1985-2000 | Boehringer Ingelheim Fonds (BIF) | Germany | fellowship | biomedical sciences |
| Brouns (2000a; NWO/ KNAW) | archive data | 1993-1994 | Dutch Organization for Scientific Research (NWO)/ Royal Dutch Academy for the Sciences (KNAW) | The Netherlands | fellowship | humanities |
| Brouns (2000b; NWO/ KNAW) | | | | | | social sciences |
| Brouns (2000c; NWO/ KNAW) | | | | | | exact sciences |
| Brouns (2000d; NWO/ KNAW) | | | | | | biology/ oceanography and earth sciences |
| Brouns (2000e; NWO/ KNAW) | | | | | | medicine |
| Dexter (2002a; Wellcome Trust) | archive | 2000-2001 | Wellcome Trust | UK | fellowship | biomedical sciences |

| Study and research funding program/ institution | Data source of the study | Cohort of application | Funding organization/ program that was examined in the study | Country location of the organization | Kind of promotion | Discipline or field |
|---|---|---|---|---|---|---|
| Dexter (2002b; Wellcome Trust) | data | | Wellcome Trust (project grant) | | grant | |
| Dexter (2002c; Wellcome Trust) | | | Wellcome Trust (program grant) | | | |
| Emery (1992; NIH) | archive data | 1990 | National Institutes of Health (NIH) | USA | fellowship | biomedical sciences |
| Friesen (1998a; MRC) | archive data | | Medical Research Council (MRC) | Canada | grant | biomedical sciences |
| Friesen (1998b; MRC) | | | Medical Research Council (MRC - Canada scholarship) | | fellowship | |
| Friesen (1998c; MRC) | | | Medical Research Council (MRC - Canada fellowship) | | | |
| Goldsmith (2002a; NSF) | archive data | 1979 | National Science Foundation (NSF) | USA | fellowship | across disciplines |
| Goldsmith (2002b; NSF) | | 1993 | | | | |
| Grant (1997a; Wellcome Trust) | archive data | 1996 | Wellcome Trust | UK | grant | biomedical sciences |
| Grant (1997b; Wellcome Trust) | | 1994-1997 | | | fellowship | |
| Grant (1997c; MRC) | | 1996 | Medical Research Council (MRC) | | grant | |
| Grant (1997d; MRC) | | 1993-1996 | | | fellowship | |
| Jayasinghe (2001; ARC) | archive data | 1996 | Australian Research Council (ARC) | Australia | grant | across disciplines |
| National Science Foundation (2005a) | archive data | 1997 | National Science Foundation (NSF) | USA | grant | across disciplines |
| National Science Foundation (2005b) | | 1998 | | | | |
| National Science Foundation (2005c) | | 1999 | | | | |
| National Science Foundation (2005d) | | 2000 | | | | |
| National Science Foundation (2005e) | | 2001 | | | | |
| National Science Foundation (2005f) | | 2002 | | | | |
| National Science Foundation (2005g) | | 2003 | | | | |
| National Science Foundation (2005h) | | 2004 | | | | |
| Over (1996; ARC) | survey | 1993 | Australian Research Council (ARC) | Australia | grant | |

| Study and research funding program/ institution | Data source of the study | Cohort of application | Funding organization/ program that was examined in the study | Country location of the organization | Kind of promotion | Discipline or field |
|---|---|---|---|---|---|---|
| Sigelman (1987; NSF) | archive data | 1985-1986 | National Science Foundation (NSF) | USA | grant | political sciences |
| Taplick (2005a; EMBO) | archive data | 2001-2004 | European Molecular Biology Organization (EMBO short-term fellowship) | Germany | fellowship | biomedical sciences |
| Taplick (2005b; EMBO) | | 1996-2004 | European Molecular Biology Organization (EMBO long-term fellowship) | | | |
| Taplick (2005c; EMBO) | | 2000-2004 | European Molecular Biology Organization (EMBO young investigator program) | | | |
| Viner (2004a; EPSRC) | archive data | 1991-1997 | Engineering and Physical Sciences Research Council (EPSRC) | UK | grant | engineering and physical sciences |
| Viner (2004b; EPSRC) | | 1995-2001 | | | | |
| Ward (1998; NHMRC) | archive data | 1994-1997 | National Health and Medical Research Council (NHMRC) | Australia | fellowship | health and medical sciences |
| Wellcome Trust (1997) | archive data | 1994-1996 | Wellcome Trust | UK | grant | biomedical sciences |
| Wenneras (1997; MRC) | archive data | 1995 | Medical Research Council (MRC) | Sweden | fellowship | biomedical sciences |
| Willems (2001a; DFG) | archive data | 1997 | German Research Foundation (DFG) | Germany | grant | biology and medical sciences |
| Willems (2001b; DFG) | | 1998 | | | | |
| Willems (2001c; DFG) | | 1999 | | | | |
| Willems (2001d; DFG) | | 2000 | | | | |
| Wood (1997a; ARC) | archive data | 1995 | Australian Research Council (ARC) | Australia | grant | biological sciences: molecular |
| Wood (1997b; ARC) | | | | | | biological sciences: plant/ animal |
| Wood (1997c; ARC) | | | | | | chemical sciences |

| Study and research funding program/ institution | Data source of the study | Cohort of application | Funding organization/ program that was examined in the study | Country location of the organization | Kind of promotion | Discipline or field |
|---|---|---|---|---|---|---|
| Wood (1997d; ARC) | | | | | | earth sciences |
| Wood (1997e; ARC) | | | | | | engineering applied sciences |
| Wood (1997f; ARC) | | | | | | humanities |
| Wood (1997g; ARC) | | | | | | physical sciences |
| Wood (1997h; ARC) | | | | | | social sciences |
| Wood (1997i; ARC) | | | | | | engineering applied sciences II |

Table 3. Studies included in the meta-analysis, with references

| Study and research funding program/ institution | Reference |
|---|---|
| Ackers (2000a-g; Marie Curie) | Ackers, L., Millard, D., Perista, H., Baptista, I., Gustafsson, V., Blomqvist, M., et al. (2000). *The participation of women researchers in the TMR Marie Curie Fellowships*. Leeds, UK: Centre for the Study of Law in Europe, University of Leeds. |
| Allmendinger (2002a-g; DFG) | Allmendinger, J., & Hinz, T. (2002). Programmed (in-)equality? Gender-specific funding of research grant proposals. *Zeitschrift für Soziologie, 31*(4), 275-293. |
| Bazeley (1998; ARC) | Bazeley, P. (1998). Peer review and panel decisions in the assessment of Australian Research Council project grant applicants: what counts in a highly competitive context? *Higher Education, 35*(4), 435-452. |
| Bornmann (2005; BIF) | Bornmann, L., & Daniel, H.-D. (2005). Selection of research fellowship recipients by committee peer review. Analysis of reliability, fairness and predictive validity of Board of Trustees' decisions. *Scientometrics, 63*(2), 297-320. |
| Brouns (2000a-e; NWO/ KNAW) | Brouns, M. (2000). The gendered nature of assessment procedures in scientific research funding: the Dutch case. *Higher Education in Europe, 25*, 193-199. |
| Dexter (2002a-c; Wellcome Trust) | Dexter, T. M. (2002). *Report on women in science, engineering and technology for Patricia Hewitt*. Retrieved June 7, 2005, from http://www.wellcome.ac.uk/assets/wtd002794.pdf. |
| Emery (1992; NIH) | Emery, J. A., Meyers, H. W., & Hunter, D. E. (1992). NIH FIRST Awards: testing background factors for funding against peer review. *Journal of the Society of Research Administrators, 24*(2), 7-28. |
| Friesen (1998a-c; MRC) | Friesen, H. G. (1998). Equal opportunities in Canada. *Nature, 391*(6665), 326-326. |
| Goldsmith (2002a-b; NSF) | Goldsmith, S. S., Presley, J. B., & Cooley, E. A. (2002). *National science foundation graduate research fellowship program. Final evaluation report*. Arlington, VA, USA: National Science Foundation (NSF). |
| Grant (1997a-d; MRC) | Grant, J., Burden, S., & Breen, G. (1997). No evidence of sexism in peer review. *Nature, 390*(6659), 438-438. |

| Study and research funding program/ institution | Reference |
|---|---|
| Jayasinghe (2001; ARC) | Jayasinghe, U. W., Marsh, H. W., & Bond, N. (2001). Peer review in the funding of research in higher education: the Australian experience. *Educational Evaluation and Policy Analysis, 23*(4), 343-346. |
| National Science Foundation (2005a-h) | National Science Foundation. (2005). *Report to the National Science Board on the National Science Foundation's Merit Review Process. Fiscal Year 2004.* Arlington, Virginia, USA: National Science Foundation (NSF). |
| Over (1996; ARC) | Over, R. (1996). Perceptions of the Australian Research Council large grants scheme: differences between successful and unsuccessful applicants. *Australian Educational Researcher, 23*(2), 17-36. |
| Sigelman (1987; NSF) | Sigelman, L., & Scioli, F. P. (1987). Retreading familiar terrain. Bias, peer review, and the NSF Political Science Program. *P S, 20*(1), 62-69. |
| Taplick (2005a-c; EMBO) | Taplick, J. (2005). Participation of women in EMBO activities. Retrieved June 7, 2006, from http://www.embo.org/gender/gender_satistics04.pdf (complete report to be published soon) |
| Viner (2004a-b; EPSRC) | Viner, N., Powell, P., & Green, R. (2004). Institutionalized biases in the award of research grants: a preliminary analysis revisiting the principle of accumulative advantage. *Research Policy, 33*(3), 443-454. |
| Ward (1998; NHMRC) | Ward, J. E., & Donnelly, N. (1998). Is there gender bias in research fellowships awarded by the NHMRC? *Medical Journal of Australia, 169*(11-12), 623-624. |
| Wellcome Trust (1997) | Wellcome Trust. (1997). *Women and peer review. An audit of the Wellcome Trust's decision-making on grants (PRISM Report No. 8).* London, UK: Wellcome Trust. |
| Wenneras (1997; MRC) | Wennerås, C., & Wold, A. (1997). Nepotism and sexism in peer-review. *Nature, 387*(6631), 341-343. |
| Willems (2001a-d; DFG) | Presentation of Silke Willems on the EMBO special meeting "The glass ceiling for women in the life sciences". Retrieved June 7, 2005, from http://www.embo.org/gender/glass_ceiling.html |
| Wood (1997a-i; ARC) | Wood, F. Q. (1997). *The peer review process.* Canberra, Australia: Australian Government Publishing Service. |

Table 4. Maximum likelihood estimates for the generalized linear mixed meta-analysis. The estimate of $\beta_1$ in the fixed part gives the effect of the applicant's gender on approval (funding) of grant applications. The random part represents the variability of regressions on gender between procedures and between studies, as well as the heterogeneity of effects.

| Effect | Parameter | Estimate | Standard error | $z$-value |
|---|---|---|---|---|
| *Fixed part* | | | | |
| Intercept | $\beta_0$ | -0.84 | 0.02 | -42.54* |
| Gender (men = -0.5; women = +0.5) | $\beta_1$ | -0.07 | 0.03 | -2.59* |
| *Random part* | | | | |
| *Peer review procedure level* | | | | |
| Intercept | var($\zeta_{0jk}$) | 0.22 | 0.03 | 14.67* |
| Gender (men = -0.5; women = +0.5) | var($\zeta_{1jk}$) | 0.04 | 0.02 | 4.87* |
| Covariance | cov($\zeta_{0jk}$, $\zeta_{1jk}$) | 0.02 | 0.02 | 1.62 |
| Correlation | cor($\zeta_{0jk}$, $\zeta_{1jk}$) | 0.24 | | |
| *Study level* | | | | |
| Intercept | var($\zeta_{0k}$) | 0.27 | 0.04 | 12.89* |
| Gender (men = -0.5; women = +0.5) | var($\zeta_{1k}$) | 0.07 | 0.02 | 5.93* |
| Covariance | cov($\zeta_{0k}$, $\zeta_{1k}$) | 0.01 | 0.02 | 0.55 |
| Correlation | cor($\zeta_{0k}$, $\zeta_{1k}$) | 0.07 | | |

* $p < .05$

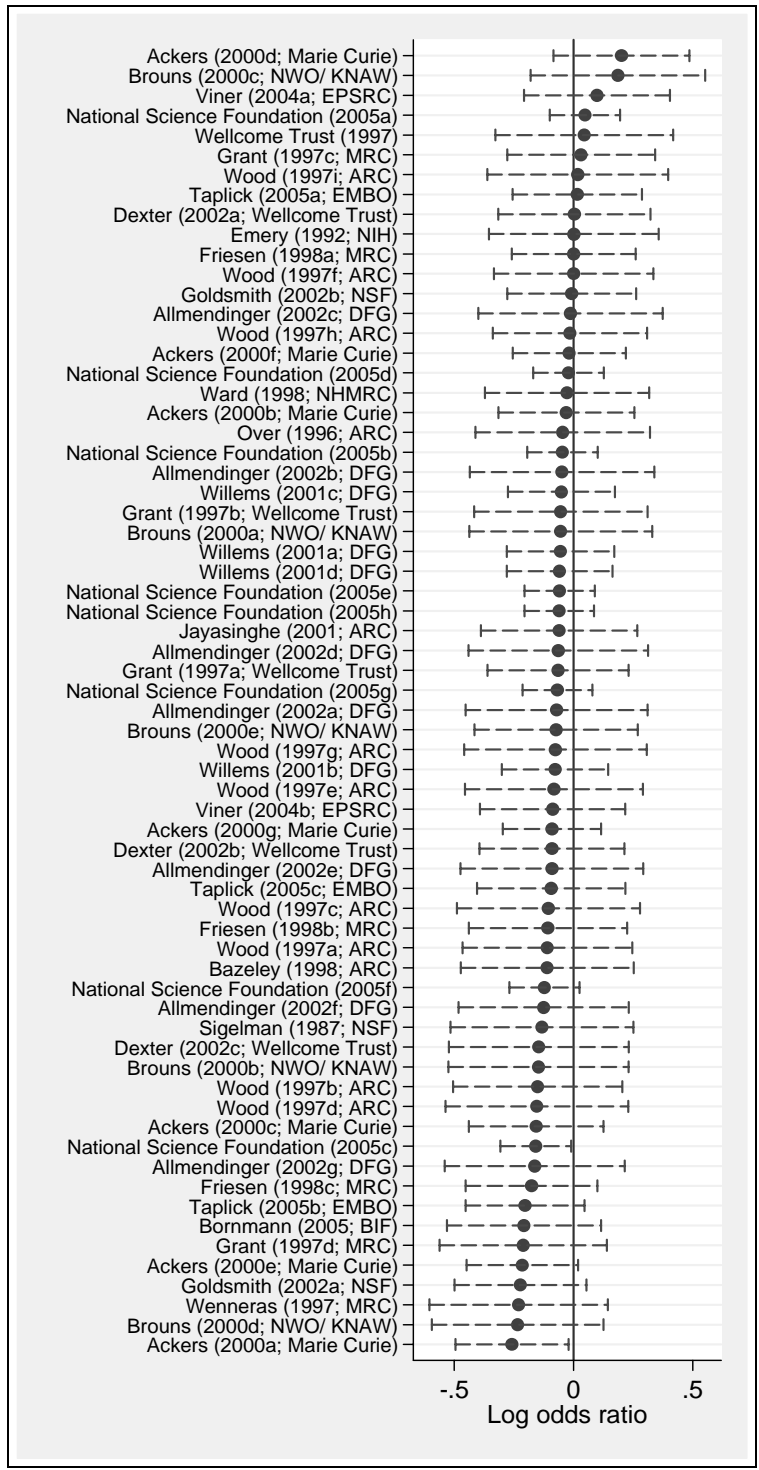

Figure 1. Mean log-odds ratios and credible intervals for each peer review procedure. Negative mean log-odds ratios point to a preference for men and positive mean log-odds ratios to a preference for women in the grant peer review.

# Acknowledgement

We would like to thank Ronald N. Kostoff, who manages science and technology evaluation at the Office of Naval Research (Arlington, VA, USA), for his support in preparing the bibliogram. We would like to thank Daryl Chubin, director of the Center for Advancing Science & Engineering Capacity (American Association for the Advancement of Science, Washington, DC), for his helpful comments on our proposal to prepare a meta-analysis on gender differences in peer review. We would like to thank Georg Matt, professor of psychology at the Medical Center of the University of California (San Diego) and Gerlind Wallon, program manager at the European Molecular Biology Organization (EMBO, Heidelberg, Germany), for their critical reading of the manuscript. The authors wish to express their gratitude to anonymous reviewers for their helpful comments.